\documentclass[a4paper,12pt]{exam}
\printanswers 
\qformat{\textbf{Problema \thequestion}\hfill}
\usepackage{color} 
\definecolor{SolutionColor}{rgb}{0.8,0.9,1} 
\shadedsolutions 

\footer{}{Pag.~\thepage}{}

\usepackage{graphicx}
\usepackage{amsmath, graphicx}

\usepackage{pgf}
\usepackage{tikz,wrapfig}
\usetikzlibrary{arrows}

\usepackage[utf8]{inputenc}
\usepackage{amsfonts, amsthm, amsmath, enumerate}
\usepackage{mathrsfs}
\usepackage[pdftex,pdfauthor={Dario Trevisan}]{hyperref}
\hypersetup{colorlinks=true}

  \newtheorem{theorem}{Theorem}
  \newtheorem{proposition}[theorem]{Proposition}
  \newtheorem{lemma}[theorem]{Lemma}
  \newtheorem{corollary}[theorem]{Corollary}
\theoremstyle{definition}
  \newtheorem{definition}[theorem]{Definition}
\theoremstyle{remark}

  \newtheorem{remark}{Remark}

 \usepackage{amsmath,color,ulem}
\newcommand{\stkout}[1]{\ifmmode\text{\sout{\ensuremath{#1}}}\else\sout{#1}\fi}

\newcommand{\re}{{\rm e }}

\newcommand{\B}{\hfill $\Box$}

\begin{document}

\title{  The Beurling-Malliavin density, the P\'olya  density and their connection }


 	\author{Rita Giuliano \footnote{Dipartimento di
		 Matematica, Universit\`a di Pisa, Largo Bruno
		 Pontecorvo 5, I-56127 Pisa, Italy (email: \texttt {rita.giuliano@unipi.it})}~~Georges Grekos\footnote{Faculté des Sciences et Techniques,  Universit\'e Jean Monnet,  23  rue Dr Paul Michelon. 42023 Saint-Etienne Cedex 2, France (email:\texttt {grekos@univ-st-etienne.fr}).  }~~Ladislav Mi\v sik       \footnote{Department of Mathematics, J. Selye University Kom\'arno, Bratislavsk\'a cesta 3322. 945 01 Kom\'arno, Slovakia (email: \texttt {misikl@ujs.sk})} \footnote {The authors are grateful to the CIRM of Marseille for the hospitality offered to them  in 2019 (under the Program {\it Research in pairs} \texttt {https://conferences.cirm-math.fr/2019-calendar.html}), when this research was at the beginning.   }   }

  \maketitle
  
\begin{abstract}

\noindent In this  paper we   present a new formulation of the Beurling-Malliavin density (Proposition \ref{prop9}). Then we  consider  the upper P\'olya  density and show how its existence is connected with the concept of subadditivity; moreover, by means of some  quantities introduced for proving Proposition \ref{prop9}, a theorem is presented   that clarifies the connection between the upper P\'olya  and the Beurling-Malliavin densities. 
In the last section we discuss the classical definition of the upper P\'olya  density and we prove a result   which seems to be new.
\end{abstract}





 \section{Introduction} 
 
 \bigskip
  \noindent
 Let $ \Lambda= \big(\lambda_n\big)_{n \in \mathbb{Z}}$ be an indexed family of real numbers. In order to solve the problem of finding the radius of completeness of $ \Lambda$ (denoted by $\mathcal{R}(\Lambda)$),   A. Beurling and P. Malliavin introduced for the first time in the paper \cite{BM} the quantity $b(\Lambda)$, defined as follows: if $\big(I_n \big)_{n \in \mathbb{Z}}$ is a sequence of disjoint intervals on $\mathbb{R}$,   call it {\it short} if
$$\sum_{n \in \mathbb{Z}}\frac{|I_n|^2}{1 + {\rm dist}^2(0, I_n)}< \infty   $$ 
(where $|I_n|=$length of $I_n$ and ${\rm dist}(0, I_n)= \min_{x \in I_n} |x| $) and {\it long} otherwise; then deﬁne
$$b(\Lambda):=\sup\{d :\exists \,\,{\rm long} \,\,(I_n)_{n\in \mathbb{Z}} \,\,{\rm such\,\, that} \,\, \#(\Lambda \cap I_n) > d|I_n|, \forall n\in \mathbb{Z}\}$$
(the particular formulation used above comes from \cite{P1} and \cite{P2}). The celebrated Theorem of  \cite{BM}  states that the radius of completeness of $ \Lambda$ is connected with $b(\Lambda)$ by the formula
$$\mathcal{R}(\Lambda) = 2 \pi b(\Lambda).$$
Later  the same quantity $b(\Lambda)$ has been studied by other authors, and various equivalent formulations have been found; see  \cite{P2} and the recent \cite{GG} for exhaustive lists of references.

\bigskip
  \noindent
Another interesting type of density is   introduced in \cite{Po} with the scope of studying gaps and singularities of power series. It is usually named as \lq\lq upper P\'{o}lya   density\rq \rq. 

\bigskip
  \noindent
In the present  paper we  first analyze the BM-density according to the definition of \cite{K2} and in particular we present a further formulation of this concept (see Proposition \ref{prop9}). Then we  study  the upper P\'olya  density and show how its existence is connected with the concept of subadditivity; moreover we prove a result (Theorem \ref{Th1}) that clarifies the connection between the upper P\'olya  and the BM-densities (by means of some  quantities introduced for proving Proposition \ref{prop9}).

\bigskip
  \noindent
In the last section of the paper we prove  that the inner limit appearing in the definition of the upper P\'olya  density, which in its original definition in \cite{Po} is calculated along the reals, can actually be calculated along the sequence of integers only. This fact seems   to have not been noticed anywhere in the past.

 \bigskip
  \noindent
  The BM-density is defined in \cite{BM}  for general families  $  (\lambda_n)_{n \in \mathbb{Z}}$ of  real numbers; on the other hand, the  P\'{o}lya density concerns sequences  of real numbers  (i.e. indexed by $n \in \mathbb{N}^*$), positive and strictly increasing ($0 \leq \lambda_1 <\lambda_2 < \cdots $ ) and such that $\lim_{n \to \infty}\lambda_n= + \infty$.    Thus, in order to compare these two concepts,  in what follows   {\it we shall confine ourselves to    sequences of the second (i.e. P\'{o}lya's)  kind .    We   also emphasize that our study concerns sequences  of real numbers and not only of integers.} In the sequel, by the term {\it sequence} we always mean a sequence of P\'{o}lya's   kind, unless otherwise specified. Notice in particular that   we are not dealing with finite sequences and sequences with repetitions.

  \bigskip
  \noindent
     The {\it counting function} of $  (\lambda_n)_{n \in \mathbb{N}^*}$ is the function
     $$F_\Lambda(t)= \begin{cases}0 & t=0\\ \#\{k  \in \mathbb{N}^*: \lambda_k \leq  t\}& t>0.
     
\end{cases}      $$For $0\leq a< b$  we have clearly
$$F_\Lambda(b)- F_\Lambda(a)= \#\{k  \in \mathbb{N}^*: \lambda_k \in (a, b]  \}.$$

  \bigskip
  \noindent
  The counting function  of a sequence $  (\lambda_n)_{n \in \mathbb{N}^*}$ is non-decreasing, is identified by the sequence itself, and viceversa;  since the definitions of the P\'{o}lya  and  B-M   densities can be formulated in terms of $F_\Lambda$, in the sequel we shall adopt the point of view  of counting functions in place of that of sequences.

  \section{The Beurling-Malliavin   density}

   \noindent
  The BM-density, firstly defined in \cite{BM}, has been studied in more recent times in \cite{K2}. This book employs a  definition which is apparently   completely different from the original one and proves  the equivalence of the two concepts. In the present paper we shall follow \cite{K2}, of course  with an adaptation to our restricted framework.

  \bigskip
   \noindent
  Let   $\mathfrak{C}$  be the family of all sequences $$\mathcal{I}= \big((a_n, b_n] \big)_{ n \in \mathbb{N}^*}$$
  of intervals in $(0, + \infty)$ such that $a_n < b_n \leq a_{n+1}$ for all $n \in \mathbb{N}^*$ and
  $$\sum_{n = 1}^\infty \Big(\frac{b_n}{a_n}-1\Big)^2= +\infty.$$
  In \cite{K2} these systems of intervals are called {\it substantial} (they are nothing but the {\it long} sequences of \cite{P1} and \cite{P2}). Then put
  $$\mathfrak{R}= \Big\{R \geq 0: \exists \,\mathcal{I}= \big((a_n, b_n] \big)_{ n \in \mathbb{N}^*}\in \mathfrak{C} ,\,\frac{F_\Lambda(b_n) - F_\Lambda(a_n)}{b_n -a_n} \geq R \, \,{ \rm for\, \, each\,\, sufficiently\,\, large\, }n\Big\}.$$
  Then the {\it Beurling--Malliavin density} $b(\Lambda)$ is defined as the supremum of $\mathfrak{R}$; in formula
  $$b(\Lambda)= \sup \mathfrak{R}.$$
  \begin{remark}
  Despite the different notation, this definition coincides with the one of  \cite{P1} and \cite{P2}.
  \end{remark}
  \bigskip
  \noindent
 For every $A \subseteq [1, + \infty]$ we shall be interested in the   subset  of $\mathfrak{C}$ defined as 
 $$\mathfrak{C}_A= 
   \Big\{ \mathcal{I}= \big((a_n, b_n] \big)_{ n \in \mathbb{N}^*}\in \mathfrak{C} : \limsup_{n \to \infty}\frac{b_n}{a_n  }\in A \Big\}.$$
  Accordingly, we shall denote
  $$\mathfrak{R}_{A}= \Big\{R \geq 0: \exists \,\mathcal{I}= \big((a_n, b_n] \big)_{ n \in \mathbb{N}^*}\in \mathfrak{C}_{A} ,\,\frac{F_\Lambda(b_n) - F_\Lambda(a_n)}{b_n -a_n} \geq R \, \,{ \rm for\, \, each\,\, sufficiently\,\, large\, }n\Big\},$$
  and finally 
  \begin{equation} \label{bA}
  b_{A}(\Lambda)= \sup \mathfrak{R}_{A}.
  \end{equation}
In the case $A = \{\alpha\}$  (i.e. $A$ is a singleton) we shall simplify this set of notation to $\mathfrak{C}_{\alpha}$, $\mathfrak{R}_{\alpha}$ and $b_{\alpha}(\Lambda)$ respectively.

\smallskip
\noindent  
For $A \subseteq B \subseteq [1, + \infty]$ we point out the obvious relations  
\begin{equation} \label{14}
\mathfrak{C}_A\subseteq\mathfrak{C}_B, \qquad \mathfrak{R}_{A}\subseteq\mathfrak{R}_{B}, \qquad  b_{A}(\Lambda)\leq b_{B}(\Lambda).
\end{equation}
  Since $a_n < b_n$ for every substantial sequence of intervals, of course
$$\limsup_{n \to \infty}\frac{b_n}{a_n }\in [1, + \infty];$$
thus, according to the preceding notation, we have $b(\Lambda)= b_{[1, + \infty]}(\Lambda)$. 

\bigskip\noindent
 For every $\mathcal{I}= \big\{(a_n, b_n] \big\}_{ n \in \mathbb{N}^*}\in \mathfrak{C} $, let
 $$\ell_{\mathcal{I}}= \liminf_{n \to \infty}\frac{F_\Lambda(b_n) - F_\Lambda(a_n)}{b_n -a_n}.$$  
 
\bigskip \noindent
 Formula \eqref{10} in the following Proposition provides an alternative definition of $b_{A}(\Lambda)$, to be used in the following section.  The  result is almost self-evident, anyway we prove it in detail.

 \begin{proposition}\label{prop9}\sl  Let $A \subseteq [1, + \infty]$ be fixed.
 Then
 \begin{equation}
 \label{10}b_{A}(\Lambda)= \sup\{\ell_{\mathcal{I}},\mathcal{I} \in \mathfrak{C}_{A} \}.
 \end{equation}
 \end{proposition}
 
 \bigskip\noindent
 {\it Proof.} Denote for simplicity $\alpha_{A}(\Lambda)= \sup\{\ell_{\mathcal{I}},\mathcal{I} \in \mathfrak{C}_{A} \}. $
  
\bigskip
\noindent
(i) We prove first that
 \begin{equation}\label{9}
 \alpha_{A}(\Lambda) \leq b_{A}(\Lambda).
 \end{equation}
  By the definition of supremum, for every $\epsilon >0 $, there exists $\mathcal{I} = \big((a_n, b_n] \big)_{ n \in \mathbb{N}^*}\in \mathfrak{C}_{A}$ such that  $$\alpha_{A}(\Lambda) - \epsilon < \ell_{\mathcal{I}}.   $$
  By definition of liminf, ultimately we have 
  $$\frac{F_\Lambda(b_n) - F_\Lambda(a_n)}{b_n -a_n}\geq \alpha_{A}(\Lambda) - \epsilon.$$
  The last inequality implies that $\alpha_{A}(\Lambda) - \epsilon\in \mathfrak{R}_{A}$, so that
  $\alpha_{A}(\Lambda) - \epsilon \leq \sup\mathfrak{R}_{A}=b_{A}(\Lambda),$ and \eqref{9} follows by the arbitrariness of $\epsilon$.

  \bigskip \noindent
  (ii) For every $R \in \mathfrak{R}_{A}$, there exists $\mathcal{I}_R = \big((a_n, b_n] \big)_{ n \in \mathbb{N}^*}\in \mathfrak{C}_{A} $ such that ultimately $$\frac{F_\Lambda(b_n) - F_\Lambda(a_n)}{b_n -a_n}\geq  R.$$ By passing to the liminf in the above inequality, we get for such $\mathcal{I }_R$ that $ R \leq  \ell_{\mathcal{I}_R} $, hence $$b_{A}(\Lambda)=\sup\mathfrak{R}_{A} \leq  \sup_R \ell_{\mathcal{I}_R} \leq \sup\{\ell_{\mathcal{I}},\mathcal{I} \in \mathfrak{C}_{A}\}= \alpha_{A}(\Lambda) . $$
  
  \B

 \section{Analysis of the Beurling-Malliavin density}

  \bigskip 
  \noindent
The aim   of the present Section is to give a precise mathematical formulation  and a rigorous proof of the intuitive feelings that, in defining the BM-density, 
\begin{itemize}
\item[(i)] all the sequences of substantial intervals  $\mathcal{I}  =\big((a_n, b_n] \big)_{ n \in \mathbb{N}^*}$ with $\limsup_{n \to \infty} \frac{b_n}{a_n}> 1$ have the same  \lq\lq status\rq\rq,   so to say (see Proposition \ref{prop3});
\item[(ii)]  for identifying the value of $b(\Lambda)$ the only important sequences of substantial intervals  $\mathcal{I} = \big((a_n, b_n] \big)_{ n \in \mathbb{N}^*}$ are those with  $\lim_{n \to \infty}\frac{b_n}{a_n}= 1$ (see Theorem \ref{prop6} here below).
\end{itemize}

 \bigskip \noindent
\begin{theorem}\label{prop6}\sl   $b(\Lambda)= b_1(\Lambda)$.
\end{theorem}

\bigskip \noindent
For proving Theorem \ref{prop6} the first step is  the following Proposition, which says that in formula \eqref{10} we can restrict ourselves to  taking the supremum of $\ell_{\mathcal{I}}$ in   particular classes $\mathfrak{C}_{A}$ (i.e.  for particular sets $A$) in place of the whole class  $\mathfrak{C}$.  
  
  \begin{proposition}\label{prop3}\sl For  $k > 1$ consider the set $A_k =[1, k]$. Then, for every $k >1$,
  $$b(\Lambda)=b_{A_k}(\Lambda).  $$
  
  \end{proposition}
  
  \bigskip \noindent
  { \it Proof.} The inequality $\geq$ is obvious by the last relation in \eqref{14}. Thus it is sufficient to show that   
  
  \bigskip \noindent
  \begin{proposition}\label{prop2}\sl For fixed $k > 1$ consider the set $A_k =[1, k]$. Then 
  \begin{itemize}
  \item[(i)] for every $\mathcal{I} = \big((a_n, b_n] \big)_{ n \in \mathbb{N}^*}\in \mathfrak{C}$, there exists $  \mathcal{J_I }\in  \mathfrak{C}_{A_k}$ such that $\ell_{ \mathcal{J_I }}\geq \ell_{\mathcal{I}}$.
  Hence $ b(\Lambda)=\sup\{\ell_{\mathcal{I}},\mathcal{I} \in \mathfrak{C}  \}\leq  \sup\{\ell_  \mathcal{J_I },\mathcal{I} \in \mathfrak{C}  \}\leq  \sup\{\ell_{\mathcal{J}},\mathcal{J} \in \mathfrak{C}_{A_k} \}=b_{A_k}(\Lambda)$.
  \item[(ii)]$\mathcal{J_I }  =    \big((c_n, d_n] \big)_{ n \in \mathbb{N}^*}$  can be chosen in such a way that $\lim_{n \to \infty}\frac{d_n}{c_n}$ exists (and     belongs to $A_k$). \end{itemize}
  \end{proposition}

  \bigskip \noindent
  For the proof of Proposition \ref{prop2} we need a lemma. 
  
   \bigskip \noindent
  \begin{lemma}\label{prop4}\sl (i) Let $k > 1$. It is possible to split any interval $(a, b]$ into some number $r$ of disjoint subintervals  
 $$(a, b]= \bigcup_{i=1}^r(c_i, d_i],$$
  where $a=c_1 < d_1 =c_2 < \cdots =c_r < d_r =b$, in such a way that, for every $i = 1, \dots, r$,
$$1 < \frac{d_{i}}{c_{i}} \leq k.$$
(ii) Let the family of intervals $\big\{(c_i, d_i], i=0, \dots, r\big\}$,   covers  the interval $(a, b]$. Then there exists an integer $  j \in \{0, \dots, r\}$ such that
  $$\frac{F_\Lambda(d_j) -F_\Lambda(c_j) }{d_j - c_j}\geq \frac{F_\Lambda(b) -F_\Lambda(a) }{\sum_{j=0}^r (d_i-c_i)}.$$
  \end{lemma}

  \bigskip \noindent  
   { \it Proof of Lemma \ref{prop4}.} (i) Take $r= \lfloor \log_k \frac{b}{a}\rfloor + 1$, $\alpha = (\frac{b}{a})^\frac{1}{r}$ and,  for $i= 1, 2, \dots, r$, define
   $$c_i = \alpha^{i-1}a, \quad d_i = \alpha^{i}a,$$
   noticing that
   $$c_1 =a, \quad d_r = b, \quad\frac{d_i}{c_i}= \alpha \leq k.$$ (ii) Assume the contrary, i.e.
  $$\frac{F_\Lambda(d_i) -F_\Lambda(c_i) }{d_i - c_i}< \frac{F_\Lambda(b) -F_\Lambda(a) }{\sum_{j=0}^r (d_i-c_i)}, \qquad \forall \, i=0, 1, 2, \dots, r.$$ Then, as the family $\big\{(c_i, d_i] \big\}$ covers $(a, b]$, 
  $$F_\Lambda(b) -F_\Lambda(a)\leq  \sum_{i=1}^r(F_\Lambda(d_i) -F_\Lambda(c_i))<  \frac{F_\Lambda(b) -F_\Lambda(a) }{\sum_{j=0}^r (d_i-c_i)} \sum_{i=1}^r(d_i-c_i)=F_\Lambda(b) -F_\Lambda(a), $$
  a contradiction. \B

  \bigskip \noindent  
  { \it Proof of Proposition \ref{prop2}.} For any $n \in \mathbb{N}^*$, by Lemma \ref{prop4} (i)  it is possible to split the interval $(a_n, b_n]$ into $r_n$ subintervals $(c_{n,i},d_{n,i} ]$ such that
\begin{equation}\label{11}
1 \leq \frac{d_{n,i}}{c_{n,i}}\leq k, \qquad i = 1, \dots r_n.
\end{equation}  
  By Lemma \ref{prop4} (ii), there exists $j_n\in \{1, 2, \dots, r_n\}$ such that
 \begin{equation}\label{12}
\frac{F_\Lambda(d_{n,j_n }) -F_\Lambda(c_{n,j_n}) }{d_{n,j_n} - c_{n,j_n}}\geq \frac{F_\Lambda(b_n) -F_\Lambda(a_n) }{b_n - a_n}.
\end{equation} 
  The sequence of numbers $$\Big(\frac{d_{n,j_n }}{c_{n,j_n }}\Big)_{n \in \mathbb{N^*}}$$
  is bounded by \eqref{11}, hence, by possibly passing to a subsequence, we can assume  that it is convergent.  It is clear that the sequence of intervals
  $$\mathcal{J}=\big((c_{n,j_n },d_{n,j_n } ]\big)_{n \in \mathbb{N}^*} $$
belongs to $\mathfrak{C}_{A_k}$ and, by \eqref{12}, verifies the relation $\ell_{\mathcal{J}}\geq \ell_{\mathcal{I}}.$ The Proposition is proved.

   \B
  
  \bigskip \noindent  
 For proving Theorem \ref{prop6} we shall use a further lemma.

  \begin{lemma}\label{prop7}\sl Let $\eta > 1$. It is possible to split any interval $(a, b]$ (with $0<a <b$) into some number $r$ of disjoint subintervals  
 $$(a, b]= \bigcup_{i=1}^r(c_i, d_i],$$
  where $a=c_1 < d_1 =c_2 < \cdots =c_r < d_r =b$ in such a way that 
$$ \frac{d_{i}}{c_{i}} =\eta,\quad i = 1, \dots, r-1, \quad \frac{d_{r}}{c_{r}} \leq\eta.$$
  
  \end{lemma}
  
   \bigskip \noindent  
  {\it Proof.} Consider the (increasing divergent) sequence $(\eta^{i-1}a)_{i \in \mathbb{N}^*}$ and denote
  $$r= \max \{i\in \mathbb{N}^*:\eta^{i-1}a< b\},$$
  which means that $\eta^{r-1}a< b$ and $\eta^{r}a \geq b$. Then put
  $$c_i = \eta^{i-1}a,   \quad   d_i =  (\eta^{i}a)\wedge b, \qquad i = 1, 2, \dots, r. $$
  Then
  $$\frac{d_{i}}{c_{i}} =\frac{\eta^{i}a}{\eta^{i-1}a}=\eta,\quad i = 1, \dots, r-1; \qquad  \frac{d_{r}}{c_{r}} = \frac{( \eta^{r}a)\wedge b}{\eta^{r-1}a}\leq \frac{ \eta^{r}a}{\eta^{r-1}a}=\eta.$$
  
 \B

   \bigskip \noindent  
 {\it Proof of Theorem \ref{prop6}}. The inequality $\geq$ is obvious by the last relation in \eqref{14}.   Thus, by Proposition \ref{prop3}, it suffices to prove that for each sequence of intervals  $\mathcal{I} = \big((a_n, b_n] \big)_{ n \in \mathbb{N}^*}\in \mathfrak{C}_{A_2}$ there exists a sequence of intervals $\mathcal{J_I}\in \mathfrak{C}_1$ such that $\ell_\mathcal{J_I} \geq \ell_\mathcal{I} $. This implies
 $$b(\Lambda)=\sup\{\ell_{\mathcal{I}},\mathcal{I} \in \mathfrak{C}  \}\leq  \sup\{\ell_  \mathcal{J_I },\mathcal{I} \in \mathfrak{C}  \}\leq  \sup\{\ell_{\mathcal{J}},\mathcal{J} \in \mathfrak{C}_{1} \}=b_{1}(\Lambda).$$
 
 \smallskip
  \noindent
  Let $n\in \mathbb{N}^*$ be fixed and put $\eta_n = (\frac{b_n}{a_n})^\frac{1}{n}$. Applying Lemma \ref{prop7} and afterwards Lemma \ref{prop4} (ii), we can construct an interval $$J_n= (c_n, d_n]\subseteq (a_n, b_n] $$
 such that
 \begin{equation}\label{13}
 1 \leq \frac{d_n}{c_n}\leq \eta_n
 \end{equation}and 
 $$ \frac{F_\Lambda(d_n) -F_\Lambda(c_n) }{d_n- c_n}\geq \frac{F_\Lambda(b_n) -F_\Lambda(a_n) }{b_n - a_n}.$$
Now notice that, since $\limsup_{n \to \infty}\frac{b_n}{a_n}\leq 2$, the sequence $(\frac{b_n}{a_n})_{ n \in \mathbb{N}^*}$ is bounded by some constant $c$, hence $\lim_{n \to \infty }\eta_n =1$ by the inequalities
$$1 \leq \liminf_{n \to \infty }\eta_n \leq \limsup_{n \to \infty }\eta_n\leq \lim_{n \to \infty }\sqrt[n]{c}=1.$$
 Thus
 $$\lim_{n \to \infty }\frac{d_n}{c_n}  =1$$
 by \eqref{13}, which concludes the proof. 
 
 \B

 \bigskip \noindent
   Let $A$ be  a subset of $[1, + \infty]$ containing 1. Then, by the last relation in \eqref{14} and by Theorem  \ref{prop6},  $$b_1(\Lambda) \leq  b_A(\Lambda)  \leq b(\Lambda)=b_1(\Lambda). $$  This proves
\begin{corollary} \label{prop8}\sl $b_{A}(\Lambda)=b(\Lambda) $ if $1 \in A$.\end{corollary}

  \section{The  upper  P\'{o}lya density $\overline{p}(\Lambda)$} 
   
  The {\it upper}  {\it P\'{o}lya density} of $\Lambda$ is the number 
  $$\overline{p}(\Lambda)= \lim_{\xi \to 1 ^-} \limsup_{x \to \infty} \frac{F_\Lambda(x) - F_\Lambda(\xi x) }{x - x \xi}.$$
In \cite{Po} it is proved that the above limit exists; see also \cite{K1}. In this section we give a new proof of this fact (see Proposition \ref{prop1}). The aim is to show how the existence of $\overline{p}(\Lambda)$ is connected with the concept of subadditivity.

  \bigskip \noindent
  The following result is a generalization of a famous lemma due to M. Fekete (see \cite{F}).
  \begin{lemma} \sl Let $g: \mathbb{R}^+\to \mathbb{R}^+$ be a function  such that there exists a continuous non--decreasing function $\phi: \mathbb{R}^+\to \mathbb{R}^+$ with the following property: for every $a, b \geq 0$,
  \begin{equation} \label{31}
  g(a+b) \leq \frac{\phi(a)}{\phi(a+b)} g(a) + \Big(1-  \frac{\phi(a)}{\phi(a+b)}\Big) g(b).
  \end{equation}
  Then $\lim_{x \to + \infty}g(x)$ exists and
  $$\lim_{x \to + \infty}g(x)= \sup_{x > 0}g(x).$$
  
  \end{lemma}
  
  \bigskip \noindent
  {\it Proof.} First it is easily seen by a simple recursive argument that, for every $a \geq 0$ and for every $n \in \mathbb{N}^*$,
  \begin{equation}\label{1}
  g(na) \leq g(a).
\end{equation}

\smallskip \noindent
   Let $y>0$ be fixed. For every   $x \in (0,y)$, let $n(x) = \lfloor \frac{y}{x}\rfloor -1$, put
  \begin{equation}\label{3}z(x) = y - n(x) \cdot x
  \end{equation}  
  and observe that, by the inequalities $a-1 <\lfloor a \rfloor \leq a$, we have
  \begin{equation}\label{2}
  x \leq z(x)\leq 2x.
  \end{equation}
  Since $y = n(x) \cdot x  + z(x)$, we get from \eqref{31} and \eqref{1} that
  \begin{align}& \label{4}
  g(y) \leq  \frac{\phi(n(x) \cdot  x)}{\phi(y)}g(x)+\Big(1-\frac{\phi(n(x) \cdot x)}{\phi(y)}\Big)g\big(z(x)\big).
  \end{align}
  
  \smallskip \noindent
   Now we distinguish the two cases (i) $\sup_{x > 0}g(x)\in \mathbb{R}$;  (ii) $\sup_{x > 0}g(x)= + \infty.$

  \bigskip \noindent
(i) It follows from \eqref{4} that
$$g(y)\leq \frac{\phi(n(x) \cdot x) }{\phi(y)}g ( x )+\Big(1-\frac{\phi(n(x) \cdot x)}{\phi(y)}\Big) \big(\sup_{x > 0}g(x)\big).$$
Let $x \to 0^+$; we have from \eqref{2} that $z(x) \to 0$, hence $n(x) \cdot x \to y$ by \eqref{3}; thus, by the continuity of $\phi$,
$$\lim_{x \to 0^+}\frac{\phi(n(x) \cdot x)}{\phi(y)} = 1$$
and
$$g(y) \leq \big(\liminf_{x \to 0^+} g(x)\big) \lim_{x \to 0^+}\frac{\phi(n(x) \cdot x)}{\phi(y)}  + \big( \sup_{x > 0}g(x)\big) \lim_{x \to 0^+}\Big(1-\frac{\phi(n(x) \cdot x)}{\phi(y)}\Big) =\liminf_{x \to 0^+} g(x).$$  
 Now we pass to the supremum in $y$ and get
 $$\sup_{y >0}g(y) \leq \liminf_{x \to 0^+} g(x) \leq \limsup_{x \to 0^+} g(x)\leq \sup_{y >0}g(y),$$
  and we are done.

   \bigskip \noindent
  (ii) We get from \eqref{4}
  $$g(y) -\Big(1-\frac{\phi(n(x) \cdot x)}{\phi(y)}\Big) g\big(z(x)\big)\leq   \frac{\phi(n(x) \cdot x)}{\phi(y)}g( x).$$
 Assume that $\liminf_{x \to 0^+}g(x)< +\infty$. Then, passing to the  $\liminf_{x \to 0^+}$, the above relation gives
 $$g(y)\leq \liminf_{x \to 0^+}g(x)< +\infty,$$ 
 hence the absurdum by passing to the sup in $y$.
 
 \B
 
 \begin{proposition} \label{prop1}\sl The limit
 $$\overline{p}(\Lambda)= \lim_{\xi \to 1 ^-} \limsup_{x \to \infty} \frac{F_\Lambda(x) - F_\Lambda(\xi x) }{x - x \xi} $$
 exists. Moreover
 \begin{equation}\label{5}
  \lim_{\xi \to 1 ^-} \limsup_{x \to \infty} \frac{F_\Lambda(x) - F_\Lambda(\xi x) }{x - x \xi}= \sup_{\xi < 1} \limsup_{x \to \infty} \frac{F_\Lambda(x) - F_\Lambda(\xi x) }{x - x \xi}.
 \end{equation} 
  \end{proposition}
 
 \bigskip \noindent
  {\it Proof.}
 The aim is to apply the previous Lemma to the function
$$g(x)=  \limsup_{y \to \infty} \frac{F_\Lambda(y) - F_\Lambda(y\re^{-x} ) }{y - y \re^{-x}}, \quad x >0,$$
 obtaining
 \begin{align*}&
 \lim_{\xi \to 1 ^-} \limsup_{y \to \infty} \frac{F_\Lambda(y) - F_\Lambda(\xi y) }{y - y \xi}= \lim_{x \to 0 ^+}\limsup_{y \to \infty} \frac{F_\Lambda(y) - F_\Lambda(y\re^{-x} ) }{y - y \re^{-x}}\\&=\sup_{x >0 }\limsup_{y \to \infty} \frac{F_\Lambda(y) - F_\Lambda(y\re^{-x} ) }{y - y \re^{-x}}= \sup_{\xi <1}\limsup_{y \to \infty} \frac{F_\Lambda(y) - F_\Lambda(\xi y) }{y - y \xi}.
 \end{align*}
 So we prove the subadditivity of $g$. We have
 \begin{align*}& 
 \frac{F_\Lambda(y) - F_\Lambda(y\re^{-(a+b)} ) }{y - y \re^{-(a+b)}}\\&= \frac{F_\Lambda(y) - F_\Lambda(y\re^{-a} ) }{y - y \re^{-a}} \cdot \frac{1 -  \re^{-a}}{1 -  \re^{-(a+b)}}+ \frac{F_\Lambda(y\re^{-a} )-F_\Lambda(y\re^{-(a+b)} )}{y\re^{-a}(1 -  \re^{-b})}\cdot \Big(1-\frac{1 -  \re^{-a}}{1 -  \re^{-(a+b)}}\Big);
 \end{align*}
 passing to the limsup in $y$ we obtain 
 $$g(a+b)\leq \frac{1 -  \re^{-a}}{1 -  \re^{-(a+b)}}g(a) + \Big(1-\frac{1 -  \re^{-a}}{1 -  \re^{-(a+b)}}\Big)g(b).$$
 Since the function $\phi(x) = 1 - \re^{-x}$ is trivially non--decreasing and continuous, the proof is concluded.
 \B

 \section{Comparison between $\overline{p}(\Lambda) $ and $ b_{ (1, + \infty]}(\Lambda)$}
  In this section we are concerned with $\mathfrak{C}_{(1, + \infty]}$, $\mathfrak{R}_{(1, + \infty]}$ and $b_{(1, + \infty]}(\Lambda)$, which we shall denote as $\mathfrak{C}_{>1}$, $\mathfrak{R}_{>1}$ and $b_{>1}(\Lambda)$ respectively for easier writing.

   \bigskip \noindent It is known that in general $\overline{p}(\Lambda)\leq  b (\Lambda)$ (see \cite{K2}). The aim of this section is to prove a more precise relation. In fact
 \begin{theorem} \label{Th1} \sl We have
  $$\overline{p}(\Lambda)=   b_{ > 1}(\Lambda).$$\end{theorem}
  
 \bigskip \noindent
We split the proof into two parts, namely Propositions \ref{prop10} and \ref{prop11}.  
 
 \begin{proposition} \label{prop10}\sl The following inequality holds true:
 $$\overline{p}(\Lambda)\leq b_{ > 1}(\Lambda).$$
  
 \end{proposition}

 \bigskip \noindent
 {\it Proof.} It suffices to show that, for every $R< \overline{p}(\Lambda)$ we have $R \leq b_{ > 1}(\Lambda).$ By definition of $ \overline{p}(\Lambda)$, there exists $\xi < 1$ such that 
 $$ \limsup_{x \to  \infty}\frac{F_\Lambda( x) - F_\Lambda(x \xi)}{ x -x \xi} > R.$$
 By definition of limsup, there exists  a  sequence   $(x_n)_{ n \in \mathbb{N}^*}$ such that
 \begin{equation}\label{7}
 \lim_{n \to \infty }{x_n}= + \infty, \quad \lim_{n \to \infty}\frac{F_\Lambda( x_n) - F_\Lambda(x_n \xi)}{ x_n -x_n \xi}=\limsup_{x \to  \infty}\frac{F_\Lambda( x) - F_\Lambda(x \xi)}{x -x \xi}> R.
 \end{equation}
 Set $n_1=1$, $y_1 = x_{n_1}= x_1$ and $n_2 = \min\{n > 1: x_n \geq \frac{x_1}{\xi}\}$ ($n_2$ exists since otherwise the sequence  $(x_n)_{ n \in \mathbb{N}^*}$ would be bounded). Put $y_2 = x_{n_2}$; then $y_2 \geq \frac{y_1}{\xi}$.
 
 \smallskip \noindent
  Assume we have constructed $n_2, \dots, n_r$ and $y_2, \dots , y_r$ such that $y_k \geq \frac{y_{k-1}}{\xi}$ for each $k = 2, \dots ,r$. Let $n_{r+1} = \min\{n >  n_r: x_n \geq \frac{x_{n_r}}{\xi}\}$ ($n_{r+1}$ exists for the same reason as above) and let $y_{r+1}=   x_{n_{r+1}} $. By this recursive construction we obtain a subsequence $(y_n)_{ n \in \mathbb{N}^*}$ of  $(x_n)_{ n \in \mathbb{N}^*}$ with the property that
  \begin{equation}\label{6}
  y_{n+1}\geq \frac{y_n}{\xi}, \qquad n=1,  2, \dots;
  \end{equation}
 now set $a_n = y_n \xi$ and $b_n = y_n$. 
 
 \noindent
 It is easy to see that the sequence of intervals  $\mathcal{I}= \big((a_n, b_n]\big)_{n \in \mathbb{N}^*}$ belongs to $\mathfrak{C}_{> 1}$:
 \begin{itemize}
  \item[(i)]$a_n < b_n \leq a_{n+1}$ since this means $ y_n \xi< y_n \leq  y_{n+1}\xi$, which is true by \eqref{6};
  \item[(ii)]$$ \sum_{n = 1}^\infty \Big(\frac{b_n}{a_n}-1\Big)^2= \sum_{n = 1}^\infty  \Big(\frac{1}{\xi}-1\Big)^2=\infty;$$\item[(iii)]$$ \limsup_{n \to \infty} \frac{b_n}{a_n }= \frac{1}{\xi}> 1.$$

  \end{itemize}
   Since ultimately
  $$\frac{F_\Lambda (b_n)-F_\Lambda (a_n)}{b_n - a_n}> R$$
  by \eqref{7}, we deduce that $R \in \mathfrak{R}_{> 1}$, hence $R \leq \sup\mathfrak{R}_{> 1}=b_{ > 1}(\Lambda). $
  
  \B
  
  \bigskip \noindent
  Now we are concerned with the reverse inequality.

 \begin{proposition}  \label{prop11}\sl We have
 $$ b_{ >1}(\Lambda)\leq \overline{p}(\Lambda).$$
 
 \end{proposition}
 
 \bigskip \noindent
 {\it Proof.} Let $R \in\mathfrak{R}_{> 1} $ and let $\mathcal{I}= \big((a_n, b_n] \big)_{ n \in \mathbb{N}^*}\in \mathfrak{C}_{> 1}$ with $\frac{\Lambda(b_n) - \Lambda(a_n)}{b_n -a_n} \geq R$. Denote $$  \limsup_{n \geq 1} \frac{b_n}{a_n}= L\geq 1.$$
 Fix $\epsilon > 0$; we have ultimately
 $$a_n > \frac{b_n}{L+\epsilon},$$
 hence
 \begin{align*}&
 R \leq \frac{F_\Lambda(b_n) - F_\Lambda(a_n)}{b_n -a_n}\leq   \frac{F_\Lambda(b_n) - F_\Lambda(\frac{b_n}{L+\epsilon})}{b_n -\frac{b_n}{L+\epsilon}}\cdot \frac{b_n (1-\frac{1}{L+\epsilon})}{b_n -a_n}= \frac{F_\Lambda(b_n) - F_\Lambda(\frac{b_n}{L+\epsilon})}{b_n -\frac{b_n}{L+\epsilon}}\cdot \frac{\frac{b_n}{a_n}(1-\frac{1}{L+\epsilon})}{\frac{b_n}{a_n} -1  } \\& \leq \frac{F_\Lambda(b_n) - F_\Lambda(\frac{b_n}{L+\epsilon})}{b_n -\frac{b_n}{L+\epsilon}}\cdot \frac{ L+\epsilon-1}{\frac{b_n}{a_n} -1  } .
 \end{align*}
 We deduce
 $$\frac{F_\Lambda(b_n) - F_\Lambda(\frac{b_n}{L+\epsilon})}{b_n -\frac{b_n}{L+\epsilon}}\geq R \cdot \frac{\frac{b_n}{a_n} -1}{L+\epsilon-1}$$
 and, passing to the limsup in $x$,
 $$\limsup_{x \to  \infty}\frac{F_\Lambda( x) - F_\Lambda(\frac{ x}{L+\epsilon})}{ x -\frac{ x}{L+\epsilon}}\geq\limsup_{x \to  \infty}\frac{F_\Lambda(b_n) - F_\Lambda(\frac{b_n}{L+\epsilon})}{b_n -\frac{b_n}{L+\epsilon}}\geq R \cdot\limsup_{n \to  \infty} \frac{\frac{b_n}{a_n} -1}{L+\epsilon-1}= R \cdot \frac{L-1}{L+ \epsilon -1}.$$
 Thus, observing that $\frac{1}{L+\epsilon} < 1$ and by Proposition \ref{prop1} (see relation \eqref{5}), we get
 $$\overline{p}(\Lambda)=\sup_{\xi < 1}\limsup_{x \to  \infty}\frac{F_\Lambda( x) - F_\Lambda(x \xi)}{ x -x \xi}\geq \limsup_{x \to  \infty}\frac{F_\Lambda( x) - F_\Lambda(\frac{ x}{L+\epsilon})}{ x -\frac{ x}{L+\epsilon}}\geq  R \cdot \frac{L-1}{L+ \epsilon -1},$$
 for every $\epsilon >0.$ Now pass to the limit as $\epsilon \to 0$ to obtain that
 $$\overline{p}(\Lambda)\geq R, \qquad \forall \,\, R \in \mathfrak{R}_{> 1},$$
 and optimizing
 $$\overline{p}(\Lambda)\geq \sup \mathfrak{R}_{> 1}=b_{ > 1}(\Lambda). $$
 
 \B 
 

\section{On the limit in the definition of the upper P\'olya  density} 
 As we have seen at the beginning of Section 4, the inner limit in the definition of $\overline{p}(\Lambda)$ is calculated as $x\to \infty$, where $x$ is a real variable. Actually  P\'olya in \cite{Po} uses the symbol $r$ (instead of $x$) without specifying where $r$ varies, but there is no reason to suppose that he didn't have the real numbers in mind. Anyway, in this Section we prove the following result:

 \begin{theorem}\label{theorem}\sl 
  
  The limit
  $$\ell:=\lim_{\eta \to 1^-} \limsup_n\frac{F(n) - F(\eta n)}{(1-\eta)n}$$
  exists and its value is $ \overline p (\Lambda).$

  \end{theorem}
  
  \noindent The proof is rather intricated and needs some preparation. In particular we need to construct a particular covering of the interval $(\xi x, x]$, i.e. a finite family of intervals  $\{(a_i, b_i]\}$ with right endpoints $b_i$ belonging to $\mathbb{N}$ and such that
  $$(\xi x, x]\subseteq \bigcup_{i}  (a_i, b_i]$$
 \bigskip \bigskip
 \subsection{Construction of an $\eta$-covering of $(\xi x, x]$}

 \bigskip \noindent 
For $x > 0$ and $\eta \in ( 0 ,1)$ let $\phi(x) = \lceil x \rceil$ and $\psi(x) = \eta x$,  $f(x) = (\phi\circ \psi)(x)= \lceil\eta  x \rceil$. 
We denote by $f^m$ the function obtained by composing $f$ with itself $m$ times, i.e.
$$f^0(x) = x; \qquad f^{m+1}(x) = (f \circ f^m )(x). $$
Similarly for $\psi^m$.

\begin{lemma}\label{Lemma2}\sl
We have the following facts:

\medskip\medskip
(i)  $t\leq  \phi(t)$, for every $t>0$;
 
 \medskip
(ii) $\psi^m(t) \leq  f^m(t)$ for every integer $m\geq 0$;

\medskip
(iii) $(f^m\circ\phi)(x)\leq \eta^mx + \sum_{k=0}^{m}\eta^k $ for every integer $m \geq 0$.\end{lemma}

 \noindent
 {\it Proof.}  (i) is evident. We prove (ii) by induction: the  case $m=0$ is obvious; the  case $m=1$ follows from (i) with $\psi(t)$ in place of $t$. For the case $m+1$  we have, 
$$\psi^{m+1}(t)= \psi(\psi^m(t))\leq \psi(f^m(t)) \leq f(f^m(t))= f^{m+1}(t)$$  
where, besides the inductive assumption, we have used the case $m=1$ and the fact that $\psi$ is increasing. 

\noindent
Now we prove (iii), again by induction. The case $m=0$ is obvious (recall that $f^0(x)=x$) and reads as $\phi(x) \leq x+1$. The inductive step uses the inductive assumption and the fact that $f$ is nondecreasing:
\begin{align*}&
f^{m+1}(\phi(x)) = f( f^m(\phi(x)))\leq f\Big(\eta^mx + \sum_{k=0}^{m}\eta^k \Big)= \Big\lceil \eta (\eta^mx + \sum_{k=0}^{m}\eta^k )\Big\rceil \\&\leq \eta \Big(\eta^mx + \sum_{k=0}^{m}\eta^k \Big)+1  = \eta^{m+1}x + \sum_{k=0}^{m+1}\eta^k.
\end{align*}
 
\hfill $\Box$

\begin{remark} By Lemma \ref{Lemma2} (iii), 
$$ \eta^mx   \leq (f^m\circ\phi)(x)\leq \eta^mx + \sum_{k=0}^{m}\eta^k $$
and, by putting $S_\eta =\sum_{k=0}^{\infty}\eta^k $, it provides the    bound
$$0 \leq \sup_{m,x}  \big((f^m\circ\phi)(x)-\eta^mx\big) \leq S_\eta.$$

\end{remark}

\bigskip
\noindent
For fixed $x$ and for every integer $i\geq 0$ put $b_i = (f^i\circ \phi)(x)$ and 
$$a_i = \eta b_i= \psi(b_i)= (\psi\circ f^i\circ \phi)(x); $$  notice that $b_i$ is an integer and $b_0 = \phi(x) = \lceil x\rceil$; notice also that, for every $i\geq 1$, 
\begin{align}& \label{2a} \nonumber
b_i = (f^i\circ \phi)(x)= f\big((f^{i-1}\circ \phi)(x)\big)= (\phi\circ \psi)\big((f^{i-1}\circ \phi)(x)\big)= \phi\big(( \psi\circ f^{i-1}\circ \phi )(x)\big)\\&= \phi(a_{i-1})=\lceil a_{i-1}\rceil.
\end{align}
Last, denote $$r(t) = \frac{1- \eta^{t+1}}{\eta^{t-1}(1-\eta)^2} , \qquad t \in [1, + \infty);$$ 
 the function $r$ is  increasing, as one can check easily by writing it in the form
 $$r(t) = \frac{1}{(1-\eta)^2}\Big(\eta e^{ t \log \frac{1}{\eta} }- \eta^2\Big).$$

\begin{lemma}\label{lemma4}\sl Let $x \geq r(q)$ for some integer $q \geq 1$. Then, for every $i = 1, \dots, q$, the following inequalities hold
  $$\eta^{i+1}x \leq a_i \leq \eta^{i }x\leq a_{i-1}\leq b_ i \leq \eta^{i-1}x\leq b_{i-1}.  $$

\end{lemma}

\noindent
 {\it Proof.}
 
 \bigskip
\noindent
 (a) That $ a_{i-1}\leq  b_i$ follows from  (ii) of Lemma \ref{Lemma2}, since, by the above definitions and putting $t = (f^{i-1}\circ \phi)(x)$, $$a_{i-1}  = \psi(t) \leq f(t)= b_i.$$

 \bigskip
\noindent
(b) Now we prove that $  \eta^{i-1}x \leq  b_{i-1}$, i.e $\psi^{i-1}(x)\leq f^{i-1}\circ \phi(x)$, which follows from   
$$\psi^{i-1}(x)\underbrace{\leq }_{Lemma \,\ref{Lemma2}, (ii)} f^{i-1}\ (x) \underbrace{\leq }_{ Lemma \,\ref{Lemma2}, (i)} f^{i-1}\circ \phi(x),$$
(since $f^m$ is non-decreasing for every $m$).

 \bigskip
 \noindent
 (c) We prove that $b_ i \leq \eta^{i-1}x$, which means $ (f^i\circ \phi)(x)\leq   \eta^{i-1}x.$

 \bigskip
 \noindent
 From Lemma \ref{Lemma2} (iii) we know that $$f^i(\phi(x))\leq \eta^ix + \sum_{k=0}^{i}\eta^k$$ and the inequality $\eta^ix + \sum_{k=0}^{i}\eta^k \leq \eta^{i-1}x$
 is equivalent to   $x \geq r(i)$; thus the claim   follows from $x \geq r(q) \geq r(i)$, recalling that $r$ is increasing.
 
 \bigskip
  \noindent
 (d) From the preceding points  (b) and (c) we have
 $$b_ i \leq \eta^{i-1}x\leq b_{i-1}$$
 and multiplying by $\eta$ we get
 $$a_i=\eta b_ i \leq \eta^{i}x\leq \eta b_{i-1}=a_{i-1}. $$
 
 \bigskip\noindent
 (e) In order to show that  $\eta^{i+1}x \leq a_i$  it suffices to multiply  by $\eta$ each side of the inequality (already proved) $\eta^i x \leq b_i$.

 \bigskip
   \noindent  
  The proof is complete.
 \hfill $\Box$
 
 \bigskip
 \begin{lemma} \label{lemma3} \sl Let $\xi\in (0,1)$ and $\eta \in (\xi, 1)$  be fixed; assume that $x \geq r(q)$ for some integer $q  \geq \frac{\log \xi}{\log \eta}$. Then  the set $\{i\geq 0: a_i > \xi x \}$ is a finite interval of integers $\{0,1, \dots, d-1\}$ for some integer $d $ with $ 1 \leq d \leq q$.

\end{lemma}
 
 \noindent
 {\it Proof.} From Lemma \ref{lemma4} we know that the sequence $\{a_i, i = 0, \dots ,q  \}$ is non-increasing. Moreover  $ a_0 = \eta \lceil x \rceil >\xi x$, while $a_q \leq \eta^q x < \xi x$, since $q \geq \frac{\log \xi}{\log \eta}$ by the assumption.  This proves the statement.
 
 \hfill $\Box$

 \bigskip
 \noindent Let $\xi\in (0,1)$ and $\eta \in (\xi , 1)$  be fixed and assume that $x \geq r(q)$ for some integer $q \geq \frac{\log \xi}{\log \eta}$. 
 From now on,    by the symbol $d$  we denote the first integer such that $a_d\leq \xi x$, i.e. $d$ verifies
 \begin{equation}\label{5a}
 a_d\leq \xi x <a_{d-1};
\end{equation}
 further   we recall that  $b_0 = \lceil x\rceil \geq x$. Thus
 $$(\xi x, x] \subseteq  \bigcup_{i=0}^d (a_i, b_i].$$
 Motivated by these remarks, we can give the following
 
 \begin{definition}   By the {\it $\eta$-covering of $(\xi x, x]$} we mean the family of intervals 
 $$\big\{(a_i, b_i], i=0, \dots, d\big\},$$  constructed as shown above.
  \end{definition}

   \bigskip \bigskip
  \subsection{Some properties of the $\eta$-covering of $(\xi x, x]$}   

   \bigskip
 \begin{lemma}\label{lemma1}\sl 
 \medskip \noindent
 Let $\xi\in (0,1)$ and $\eta \in (  { \xi} , 1)$  be fixed;    assume that   $x \geq \frac{1- \eta^2  \xi}{ \xi   (1-\eta)^2}. $   Let $$\{(a_i, b_i], i= 0, 1, \dots, d\}$$ be the $\eta$-covering of $(\xi x, x]$. Then
  $ d\leq  
  \lceil  \frac{\log \xi}{\log \eta}  \rceil. $ 
  
 
 \end{lemma}
 
   \noindent
 {\it Proof.}  In the proof we denote $\lceil  \frac{\log \xi}{\log \eta}\rceil=:q$ for simplicity.  
 
  \noindent
 Notice   that
 $$x\geq\frac{1- \eta^2  \xi}{ \xi   (1-\eta)^2}= r\Big(\frac{\log \xi}{\log \eta} +1\Big)>r(q).$$
 Hence Lemma \ref{lemma4} is in force for the integers $  1, \dots , q$. Applying this lemma with $i=q$  and by the definition of the ceiling function, we find 
 \begin{equation*}
 a_q\leq \eta^q x \leq \xi x . 
 \end{equation*}
 By the definition of $d$, this relation says that $d\le q$. 
 
  
 \hfill $\Box$

 \bigskip 
 \noindent
 \begin{remark}\label{remark}
 Actually, we can prove even more, precisely that $d \in  \{  \lfloor \frac{\log \xi}{\log \eta}  \rfloor-1, \lfloor  \frac{\log \xi}{\log \eta}  \rfloor \}$. The proof of this fact is rather complicated, and is postponed in the last Section 7. Here Lemma \ref{lemma1} will be sufficient for our scopes.
 \end{remark}

 \begin{lemma}\label{lemma5}\sl Let $\xi\in (0,1)$ and $\eta \in (\xi, 1)$  be fixed. Denote
 $$M=M(\xi, \eta):=  \frac{1- \eta^2  \xi}{ \xi    (1-\eta)^2}\vee \frac{ \frac{\log \xi}{\log \eta} +2}{\eta^2 \xi (1-\eta)}.$$
 Let $x \geq M $ and let $\{(a_i, b_i], i= 0, 1, \dots, d\}$ be the $\eta$-covering of $(\xi x, x]$).  
  
 \bigskip

 $$\sum_{i=0}^d (b_i -a_i)\leq(1-\xi \eta^3)x. $$

   \end{lemma}

  \noindent
 {\it Proof.} 
 Once more  by Lemma \ref{lemma4}, for every $i$ the two   intervals $(a_{i-1}, b_{i-1}]$ and $(a_{i}, b_{i}]$ overlap  on the interval $(a_{i-1}, b_{i}]$, the length of which  is   $b_{i}-a_{i-1}= \lceil a_{i-1}\rceil- a_{i-1}\leq 1,$ due to \eqref{2a}. Hence (recall that $b_0= \lceil x\rceil$)
$$\sum_{i=0}^d (b_i -a_i)\leq b_0 - a_d + d = \lceil x\rceil- a_d + d\leq x+1 -a_d +d\leq x+1 - \eta^{d+1 }x+d, $$ 
 where we have  used the left-hand inequality in  Lemma \ref{lemma4}. Continuing and using Lemma \ref{lemma1}, we find
  $$x+1 - \eta^{d+1 }x+d= (1-\eta^{d+1 })x + d +1 \leq (1-\eta^{2 }\xi)x + d +1\leq   (1-\eta^{2 }\xi)x +    \frac{\log \xi}{\log \eta}   +2\leq (1-\eta^{3 }\xi)x,$$
where   the first and second inequalities come  from $d  \leq \frac{\log \xi}{\log \eta}+1$ (by Lemma \ref{lemma1}),  and  the last one holds  since $x \geq  M $.

 \hfill $\Box$
 
  \bigskip \bigskip
 \subsection{ Use of the $\eta$-covering of $(\xi x, x]$}    
 Now we are in a position to prove Theorem \ref{theorem}.

 \bigskip
 \noindent
 First notice the obvious relations
  \begin{align*}&
  \liminf_{\eta \to 1^-}\Big(\limsup_n \frac{F(n) - F(\eta n)}{(1-\eta)n}\Big)\leq \limsup_{\eta \to 1^-}\Big(\limsup_n \frac{F(n) - F(\eta n)}{(1-\eta)n}\Big)\\&\leq \sup_{\eta \in (0,1)}\Big(\limsup_n \frac{F(n) - F(\eta n)}{(1-\eta)n}\Big)\leq \sup_{\xi \in (0,1)}\Big(\limsup_x  \frac{F(x)- F(\xi x)}{(1-\xi  )x}\Big);
  \end{align*}
 hence, in order to prove   Theorem \ref{theorem}, it suffices to show that

 \begin{proposition}\label{proposition8}\sl We have
 $$\liminf_{\eta \to 1^-}\Big(\limsup_{n \to \infty} \frac{F(n) - F(\eta n)}{(1-\eta)n}\Big)\geq\sup_{\xi \in (0,1)}\Big(\limsup_{x \to \infty}  \frac{F(x)- F(\xi x)}{(1-\xi  )x}\Big).$$
 
 \end{proposition}
   
 \noindent
 {\it Proof.}  We shall use the following   Lemma, the proof of which is postponed at the end.

 \begin{lemma}\label{7a}\sl  Let $\xi\in (0,1)$ and $\epsilon \in (0,1)$ be fixed. Then there exists $\delta = \delta(\xi, \epsilon)$  with the following property: for every $\eta \in  (\delta, 1)$ there exists $M=M(\xi, \epsilon, \eta)$ such that for every $x>M$ there exists  an integer $n= n(\xi, \epsilon, \eta,x)$  
 $$\frac{F(n) - F(\eta n)}{(1-\eta)n}\geq  \frac{F(x)- F(\xi x)}{(1-\xi  )x}(1-\epsilon).$$

 \end{lemma}

 \bigskip
 \noindent
  Now, in order to prove Proposition \ref{proposition8}, observe that in Lemma \ref{7a} the integer $n$ depends on $x$, so here we denote it by $n_x$. Thus, passing to the limsup as $x \to \infty$ in the relation of Lemma \ref{7} we obtain that, for any fixed $\xi$ and $\epsilon$, there exists $\delta(\xi, \epsilon)$ such that, for every $\eta \in (\delta(\xi, \epsilon), 1)$,
$$\limsup_{x \to \infty}\frac{F(n_x) - F(\eta n_x)}{(1-\eta)n_x}\geq \limsup_{x \to \infty} \frac{F(x)- F(\xi x)}{(1-\xi  )x}(1-\epsilon).$$
  Consequently we have also
  $$\limsup_{n \to \infty}\frac{F(n ) - F(\eta n )}{(1-\eta)n }\geq \limsup_{x \to \infty} \frac{F(x)- F(\xi x)}{(1-\xi  )x}(1-\epsilon).$$
  Denote provisorily
  $$A(\xi)= \limsup_{x \to \infty} \frac{F(x)- F(\xi x)}{(1-\xi  )x}; \qquad B(\eta)=\limsup_{n \to \infty}\frac{F(n ) - F(\eta n )}{(1-\eta)n }.$$
  Then we have  that, for any fixed $\xi$ and $\epsilon$, there exists $\delta(\xi, \epsilon)$ such that, for every $\eta \in (\delta(\xi, \epsilon), 1)$,
  $$B(\eta)\geq A(\xi) (1-\epsilon);$$
  hence, for every $\xi$ and $\epsilon$,
  $$ \liminf_{\eta \to 1^-} B(\eta)\geq A(\xi) (1-\epsilon),$$
  and now the statement follows by the arbitrariness of $\xi$ and $\epsilon$.
  
  \bigskip
  \noindent
  This concludes the proof of Proposition \ref{proposition8}.
  
  \hfill $\Box$
  
   \bigskip
    \noindent
  It remains to give the proof of Lemma \ref{7a}. To this extent, we need the following

     \begin{lemma}\label{6a}\sl Let $\xi\in (0,1)$ and $\eta \in (\xi, 1)$  be fixed.   Then there exists $M=M(\xi, \eta  )>0$ such that, for every  $x >M$, there exists an integer $n = n(\xi,x,\eta  )$ that verifies the inequality
$$\frac{F(n) - F(\eta n)}{(1-\eta)n}\geq \frac{F(x)- F(\xi x)}{(1-\xi \eta^3)x}.$$\end{lemma}

    \bigskip
    \noindent{\it Proof of Lemma \ref{6a}.} Let $M$   be as in  Lemma   \ref{lemma5}. If $x>M$, the $\eta$-covering of the interval $(\xi x, x]$,
$\{(a_i ,b_i], i= 0, \dots,  d\}$ is such that
$$\sum_i (b_i -a_i)\leq(1-\xi \eta^3)x. $$
By Lemma \ref{prop4}(ii), there exists $j \in\{0, \dots, d\} $ such that
$$\frac{F(b_j)-F(a_j)}{b_j - a_j } \geq\frac{F(x)- F(\xi x)}{\sum_i (b_i -a_i) } .$$
  By construction of an $\eta$-covering, $  b_j$  is an integer and $a_j = \eta b_j$; thus, if we call  $b_j = n$, we have obtained
 $$\frac{F(n) - F(\eta n)}{(1-\eta)n}\geq \frac{F(x)- F(\xi x)}{\sum_i (b_i -a_i)}\geq \frac{F(x)- F(\xi x)}{(1-\xi \eta^3)x}. $$
 
 \hfill $\Box$

  \bigskip
\noindent
 {\it Proof of Lemma \ref{7a}.} 
Since
 $$\lim_{\eta \to 1^-}\frac{1-\xi}{1-\xi \eta^3}=1,$$
   there exists $\overline\delta(\xi, \epsilon)$ such that, for every $\eta \in (\overline\delta, 1)$
 $$ 1-\xi \eta^3 < \frac{1-\xi}{1-\epsilon} .$$
Now take $\delta(\xi, \epsilon)= \overline\delta(\xi, \epsilon)\vee \xi$ and apply Lemma     \ref{6a}  to conclude.
  \hfill $\Box$
  
\section{On the number $d$}
This Section is devoted to the proof of the announced result concerning $d$ (see Remark \ref{remark}). Precisely 
  \begin{proposition}\label{prop17}\sl     
   
   (i) Assume that  $\frac{\log \xi}{\log \eta} $ is not an integer and 
  \begin{equation}\label{35}
   x \geq  \frac{1- \eta ^2 \xi}{ \xi  (1-\eta)^2}\vee \frac{ \lfloor\frac{\log \xi}{\log \eta}\rfloor +2}{\xi - \eta^{ \lfloor\frac{\log \xi}{\log \eta}\rfloor +1}}\vee   \frac{1}{  \eta^{  \lfloor\frac{\log \xi}{\log \eta}\rfloor }-\xi}.
\end{equation}   
  Then $d =   \lfloor\frac{\log \xi}{\log \eta}\rfloor  $.
  
   \bigskip 
 \noindent
   (ii) If  $\frac{\log \xi}{\log \eta} $ is an integer and $$x \geq  \frac{1- \eta^2  \xi}{ \xi    (1-\eta)^2}\vee \frac{  \frac{\log \xi}{\log \eta}  +2}{\xi(1 - \eta)},$$
   then $d \in  \{  \frac{\log \xi}{\log \eta}  -1 , \frac{\log \xi}{\log \eta} \}.$
   
   \end{proposition}
   
   \bigskip
  \noindent
 {\it Proof.} In the proof we denote $ \lfloor \frac{\log \xi}{\log \eta}\rfloor=:p$ for simplicity. Notice also that  Lemma \ref{lemma4} is applicable for $i=1, \dots p+1$, since, in both (i) and (ii), $$  x\geq   \frac{1- \eta^2  \xi}{ \xi    (1-\eta)^2}= r\Big(\frac{\log \xi}{\log \eta}+1\Big)\geq r(p+1).$$  
Last,   in case (i), $\eta^{p+1}< \xi< \eta^{p}$, hence   the denominators $\xi - \eta^{p +1}$ and $\eta^{p}-\xi$   appearing in \eqref{35} are strictly positive.

\bigskip
  \noindent
The trivial observation that $\eta^m  x\leq b_m$ and Lemma \ref{Lemma2} (iii) yield for every $m \geq 0$
 $$\eta^m x \leq b_m \leq \eta^m x + m +1,$$
 and by \eqref{2a}, for every $m \geq 0$
 $$b_m -1 < a_{m-1}\leq b_m.$$
 Putting $m=p$ and $m=p+1$ we obtain in turn
 \begin{equation}\label{3a}
 \eta^p x \leq b_p \leq \eta^p x + p +1, \qquad b_p -1 < a_{p-1}\leq b_p
 \end{equation}
 and
 \begin{equation}\label{4a}
  \eta^{p+1} x \leq b_{p+1} \leq \eta^{p+1} x + p +2, \qquad b_{p+1} -1 < a_{p}\leq b_{p+1}.
 \end{equation}
Since $x\geq  \frac{ p +2}{\xi - \eta^{p+1}}$, we have, by \eqref{4a} and Lemma \ref{lemma4},
\begin{equation}\label{36}
  a_p \leq b_{p+1} \leq \eta^{p+1} x + p +2 \leq \xi x,
 \end{equation}
 which says that $d \leq p$.  
 
 \bigskip
 \noindent
 Now we distinguish the two cases (i) and (ii).
 
 \bigskip
 \noindent
 (i) If $\frac{\log \xi}{\log \eta}   $ is not an integer, by  \eqref{3a},
 $$a_{p-1}> b_p -1 \geq \eta^p x-1 >  \xi x,$$
 Then, by \eqref{36},
 \begin{equation}\label{8a}
 a_p \leq \xi x <  a_{p-1},
 \end{equation}
 implying that $d = p$. 
 
 \bigskip
 \noindent
 (ii)  If $\frac{\log \xi}{\log \eta}   $ is  an integer, then $\xi x = \eta^p x $ and, by Lemma \ref{lemma4}, $  \eta^p x \leq a_{p-1}.$ 
 
 \bigskip
 \noindent
  If $\xi x (= \eta^p x) < a_{p-1}$, we have the same inequalities as before in \eqref{8a}
 and $d=p$ again.
 Otherwise $\eta^p x  (= \xi x ) = a_{p-1}$. We prove that in this case $d= p-1$, which means that
 \begin{equation}\label{34}
  a_{p-1}\leq \xi x < a_{p-2}. 
 \end{equation}
 The proof of \eqref{34} follows from
 \begin{lemma}\label{lemma6}\sl Assume that for some $i\geq 0$ we have $\eta^{i+1} x=a_{i} $. Then 
 $$ \eta = \frac{b_1}{b_0}, \qquad x= b_0 , \qquad a_j = \frac{b_1^{j+1}}{b_0^j}, \qquad  j= 0, \dots, i.$$
 
 \end{lemma}
 
  \noindent
 {\it Proof } (of Lemma \ref{lemma6}). First, recall the definition of every $a_j$, i.e. 
 \begin{equation}\label{8b}
  a_j = \eta b_j.
 \end{equation}
 Now, from $ a_{i} = \eta^{i+1} x  $ we deduce by \eqref{8b} that $b_i= \eta^i x $. Since $\eta^i x\leq a_{i-1} \leq b_i$ (by Lemma \ref{lemma4}), we obtain that $a_{i-1}=b_i$, whence $ \eta b_i=\eta a_{i-1}$ and (by \eqref{8b}) $a_i = \eta a_{i-1}$. Repeating this argument we get that $a_j= \eta a_{j-1}$ for every $j = 1, \dots, i$; from $  a_1 = \eta a_0$ and  \eqref{8b} (applied to $a_0$ and  $a_1$)  we deduce  that $b_1 = \eta b_0$, hence $\eta =\frac{b_1}{b_0} $, which implies that 
 $$a_0 = \eta b_0=  \frac{b_1}{b_0}\cdot b_0= b_1, \qquad a_1 = \eta a_0= \frac{b_1}{b_0}\cdot b_1= \frac{b^2_1}{b_0} , \dots  , a_j =\frac{b^{j+1}_1}{b^j_0}, \qquad j= 2, \dots, i;$$
 further
 $$x = \frac{a_i}{\eta^{i+1}}=  \frac{b^{i+1}_1}{b^i_0}\cdot \frac{b_0^{i+1}}{ b_1^{i+1}}= b_0.$$ \hfill $\Box$
 
 \bigskip
 \noindent
 We go back to the proof of Proposition \ref{prop17}; 
 the left inequality in \eqref{34} holds obviously (it is nothing but the equality  $a_{p-1}   = \eta^p x =\xi x$). Concerning the right one, by the same equality  and   Lemma \ref{lemma6},  it is equivalent to
  $\frac{b_1^{p }}{b_0^{p-1}} < \frac{b_1^{p-1}}{b_0^{p-2}}$ 
 and, after simplification, to
  $ \frac{b_1}{b_0}< 1,$ 
 which is true since  it is nothing but $\eta < 1 $ (by   Lemma \ref{lemma6} again).

 \hfill $\Box$

 \end{document}